\documentclass[12pt,leqno]{article}
\usepackage{amsmath,amsthm,amscd,amssymb}
\usepackage{latexsym}
\newtheorem{remark}{Remark}[section]
\newcommand{\bremark}{\begin{remark} \em}
\newcommand{\eremark}{\end{remark} }

\begin{document}
\parindent 15pt
\renewcommand{\theequation}{\thesection.\arabic{equation}}
\renewcommand{\baselinestretch}{1.15}
\renewcommand{\arraystretch}{1.1}
\def\disp{\displaystyle}
\title{\large Faber-Krahn inequality for Robin problem involving p-Laplacian \footnotetext{This work is supported by NNSF
of China(Grant: No.10671064)} \footnotetext{}}
\author{{\small Qiuyi Dai}\\
{\small Department of Mathematics, Hunan Normal University}\\
 {\small Changsha Hunan 410081, P.R.China}\\
{\small E-mail: daiqiuyi@yahoo.com.cn}\\
{\small Yuxia Fu}\\
{\small Department of Applied Mathematics, Hunan University}\\
 {\small Changsha Hunan 410082, P.R.China}\\
{\small E-mail: fu\_yuxia@yahoo.com.cn}\\
\date{}}
\maketitle

\abstract{\small The eigenvalue problem for the p-Laplace operator
with Robin boundary condition is considered in this paper. A
Faber-Krahn type inequality is proved. More precisely, it is shown
that amongst all the domains of fixed volume, the ball has the
smallest first eigenvalue.}



\section*{Preface}

The result of this paper was announced at first on a conference held
at the Wuhan Institute of Physics and Mathematics in May 2007. The
final version of this paper was finished in September 2008 when the
first author worked as a research fellow at The Australian National
University. As soon as we completed our paper, we sent a copy of our
preprint to D.Daners (see item 6 in the reference of \cite{sBDan})
since our result is related to a previous paper \cite{sDan} of him.
Five months later, D.Bucur and D.Daners give an alternative proof of
our result in February 2009. Though the paper has been published
(see \cite{sBDan}), their proof depends completely on Proposition
2.2, Corollary 2.3 and Proposition 2.7 of this paper which, to our
knowledge, can not been found in other materials. For completeness
of their proof and the reader's convenience, our paper will be
published here.

\section*{1 Introduction}

\setcounter{section}{1}

\setcounter{equation}{0}

\noindent Let $\Omega \subset R^N(N\geq 2)$ be an open bounded
smooth domain, we consider the following eigenvalue problem
\begin{equation}
\left\{\begin{array}{ll}\label{eq1}
-{\rm{div}}(|\nabla u|^{p-2}\nabla u)=\lambda |u|^{p-2}u & \mbox{in}\ \ \Omega,\\
|\nabla u|^{p-2}\frac{\partial u}{\partial \nu}+\beta |u|^{p-2}u=0 & \mbox{on}\ \ \partial\Omega,\\
\end{array}
\right.
\end{equation}
where $1<p<+\infty$, $\nu$ is the outward unit normal of
$\partial\Omega$ and $\beta$ is a non-negative constant.

The p-Laplacian ${\rm{div}}(|\nabla u|^{p-2}\nabla u)$ arises in
many applications such as non-Newtonian fluids, quasi-regular and
quasi-conformal mapping theory and Finsler geometry etc. An
important special case of the p-Laplacian is the well known
Laplacian $\Delta u={\rm{div}}(\nabla u)$ which corresponds to
$p=2$. Problem (\ref{eq1}) is called Dirichlet when
$\beta=+\infty$, Neumann when $\beta=0$, and Robin when
$0<\beta<+\infty$.

The main purpose of this paper is to prove a Faber-Krahn type
inequality for the Robin problem of the p-Laplacian. This
inequality says that amongst all the domains of fixed volume, the
ball has the smallest first eigenvalue. The study of this kinds of
inequalities can be traced back to 1877 \cite{sRay}. Let $B$
denote a ball in $R^N$, and $\lambda^D_1(\Omega)$ denote the first
eigenvalue of the following eigenvalue problem
\begin{equation}\label{eq2}
\left\{
\begin{array}{ll}
-\Delta\psi=\lambda\psi &x\in\Omega,\\
\psi=0 &x\in\partial\Omega.
\end{array}
\right.
\end{equation}
Rayleigh \cite{sRay} conjectured that
\begin{equation}\label{eq3}
\lambda^D_1(\Omega)\geq\lambda^D_1(B)\ \ \ \mbox{for}\
 \ \Omega\subset R^N\ \ \ \mbox{with}\ |\Omega|=|B|,
\end{equation}
and the equality hold if and only if $\Omega=B$. This conjecture
was proved independently by Faber \cite{sFab} and Krahn
\cite{sKra,sKra1} in the 1920's by making use of Schwartz
symmetrization. Since then, the inequality (\ref{eq3}) was known
as Faber-Krahn inequality. In 1999, a proof of Faber-Krahn type
inequality for the Dirichlet problem of the p-Laplacian was given
by T. Bhattacherya \cite{sBha}. Recently, Faber-krahn type
inequality was generalized to Robin problem of the Laplacian by
M.H.Bossel \cite{sBos} for dimension $N=2$, and by D.Daners
\cite{sDan} for dimension $N\geq 3$ but left the equality case
open. A little bit later, D.Daners and J.Kennedy complete the
proof of equality case in \cite{sDanK}. Note that the
generalization of the Faber-Krahn inequality from Dirichlet
problem to Robin problem is not trivial as, unlike in the
Dirichlet problem, the first eigenvalue of Robin problem is not
monotone as the domain expands (see \cite{sGio}). For more
information of the Faber-Krahn type inequality on manifold , we
refer to \cite{sIssa}.

Since the level surface of the first eigenfunction of Robin problem
intersects with the boundary $\partial\Omega$, the Schwartz
symmetrization of the first eigenfunction generally does not
decrease its Dirichlet integral and hence the Schwartz
symmetrization method does not apply to the proof of Faber-Krahn
inequality for Robin problem. Therefore, new approach must be
employed in the proof of the Faber-Krahn inequality for Robin
problem. The two crucial tools used by D.Daners \cite{sDan} to prove
the Faber-Krahn inequality for Robin problem of the Laplacian are
the Bessel functions and a new formula for the first eigenvalue by
making use of level sets of the corresponding eigenfunction. To
prove the Faber-Krahn type inequality for Robin problem of the
p-Laplacian with $p\neq 2$, we mainly face two difficulties. One is
the lack of Bessel functions and the other is the degeneracy of the
operator. The tools we use to overcome these difficulties are some
new abstract propositions of the first eigenfunction and some
approximation procedure.  The main results of this paper can be
stated as the following

\vskip 0.1in

{\em{\bf Theorem 1.1.} Let $1<p<+\infty$ and $\lambda_1(\Omega)$
be the first eigenvalue of problem (1.1) with $0<\beta<+\infty$.
If $B$ is an open ball such that $|B|=|\Omega|$, then
$\lambda_1(B)\leq\lambda_1(\Omega)$.}

\vskip 0.1in

{\bf Remark.}\ Theorem 1.1 is proved under the assumption that
$\Omega$ is smooth. However, by an approximation method similar to
that used in \cite{sDan}, we can prove that Theorem 1.1 is still
true for the domains of Lipschitz type.

\vskip 0.1in

{\em{\bf Theorem 1.2.} Let $1<p<+\infty$ and $B$ be a ball
satisfying $|B|=|\Omega|$. If $\lambda_1(\Omega)=\lambda_1(B)$,
then, up to a translation, we have $\Omega=B$.}

\vskip 0.11in

We also point out here that a symmetry result due to Gidas, Ni and
Nirenberg \cite{sGNN} plays a crucial role in the proof of Theorem
1.2 when $p=2$ (see D.Daners and J.Kennedy \cite{sDanK}). However,
this kind of result is not available for p-Laplace equation when
$p>2$ and $p\neq N$ (see however \cite{sKeP} for the case $p=N$).
Fortunately, we can prove a symmetry result needed in the proof of
Theorem 1.2 in the special case of eigenvalue problem, though we
can not prove more general symmetry result as in \cite{sGNN}.

The contents of the rest of this paper are as follows:\ $\S 2$\
The First Eigenvalue and Eigenfunction. $\S 3$\ Level Sets Formula
of $\lambda_1(\Omega)$. $\S 4$\ The lower bound of
$\lambda_1(\Omega)$. $\S 5$\ Proof of Theorem 1.1. $\S 6$\ Proof
of Theorem 1.2.

\section*{2 The First Eigenvalue and Eigenfunction}

\setcounter{section}{2}

\setcounter{equation}{0}

\noindent In this section, we give definition and some properties
of the first eigenvalue and its corresponding eigenfunction of
problem (\ref{eq1}). We focus on the case $0<\beta<+\infty$, since
the case $\beta=+\infty$ has been resolved and the case $\beta=0$
is trivial.

Let $\mathcal{K}=\{u\in W^{1,p}(\Omega);\|u\|_{L^p(\Omega)}=1\}$.
Let
\begin{equation}\label{eq4}
\lambda_1(\Omega)=\inf\{\int_{\Omega}|\nabla
u|^p+\beta\int_{\partial\Omega}|u|^p;\ u\in\mathcal{K}\}
\end{equation}
be the first eigenvalue and $\psi$ the corresponding eigenfunction
of problem (\ref{eq1}).

\vskip 0.2cm

{\em {\bf Proposition 2.1.} Let $\lambda_1(\Omega)$ be defined as in
(\ref{eq4}). Then $\lambda_1(\Omega)>0$ can be achieved by some
positive function $\psi$.}

\vskip 0.2cm

{\bf Proof.} Define functional $\Phi(u)$ on $\mathcal{K}$ by
$$\Phi(u)=\int_{\Omega}|\nabla u|^p dx +\beta\int_{\partial \Omega}|u|^p dx.$$
It is obvious that $\Phi(u)$ is a convex functional. By Theorem
1.3 in Chapter 5 of \cite{sGuoD}, $\Phi(u)$ is weakly lower
semi-continuous on $\mathcal{K}$. Let $\{u_j\}_{j=1}^{\infty}$ be
a minimum sequence of $\lambda_1$ on $\mathcal{K}$, that is,
$\int_{\Omega}|u_j|^p dx=1$ and
$$\int_{\Omega}|\nabla u_j|^p dx
+\beta\int_{\partial \Omega}|u_j|^p dx\rightarrow
\lambda_1(\Omega),\ \ \rm{as}\ \ j\rightarrow +\infty.$$

Since $\{u_j\}$ is bounded in $W^{1,p}(\Omega)$ and the embedding
$W^{1,p}(\Omega)\hookrightarrow L^p(\Omega)$ is compact, there
exists $u\in W^{1,p}(\Omega)$ such that
$$u_j\rightharpoonup u \ \ \rm{weakly\ in}\ W^{1,p}(\Omega),$$
$$u_j\rightharpoonup u \ \ \rm{strongly\ in}\ L^p(\Omega).$$

Hence, by the weakly lower semi-continuity of $\Phi(u)$, we have
$$\Phi(u)\leq \underline{\lim}_{j\rightarrow+\infty} \Phi(u_j)=\lambda_1(\Omega).$$
On the other hand, we have $\lambda_1(\Omega)\leq\Phi(u)$ due to
$u\in\mathcal{K}$ and the definition of $\lambda_1(\Omega)$. Thus
$\lambda_1(\Omega)=\Phi(u)$. Let $\psi=|u|$ , it is easy to check
that $\lambda_1(\Omega)=\Phi(\psi)$. Moreover, $\psi$ is positive
in $\bar{\Omega}$ by the strong maximum principle (see Lemma 2.6
below) . Thus, we complete the proof of proposition 2.1.

\vskip 0.2cm

{\em {\bf Proposition 2.2.}\ Let $\lambda_1(\Omega)$ be the first
eigenvalue of problem (\ref{eq1}). Then $\lambda_1(\Omega)$ is
simple in the sense that if $\psi_1>0$ and $\psi_2>0$ are two
eigenfunctions corresponding to $\lambda_1(\Omega)$, then
$\psi_2=C\psi_1$ and $C$ is a constant.}

{\bf Proof.} Suppose that $\psi_1$ and $\psi_2$ are two
eigenfunctions corresponding to $\lambda_1(\Omega)$, and $\psi_1,\
\psi_2>0$. Then $\psi_i,\, i=1,2$ satisfy
\begin{equation}\label{eq5}
\left\{\begin{array}{ll}
-{\rm{div}}(|\nabla\psi_i|^{p-2}\nabla\psi_i)=\lambda |\psi_i|^{p-2}\psi_i & \mbox{in}\ \Omega,\\
|\nabla\psi_i|^{p-2}\frac{\partial\psi_i}{\partial \nu}+\beta |\psi_i|^{p-2}\psi_i=0 & \mbox{on}\ \partial\Omega.\\
\end{array}
\right.
\end{equation}

Let
$$\eta_1=\psi_1-\psi_2^p\psi_1^{1-p}=\frac{\psi_1^p-\psi_2^p}{\psi_1^{p-1}},
\ \
\eta_2=\psi_2-\psi_1^p\psi_2^{1-p}=\frac{\psi_2^p-\psi_1^p}{\psi_2^{p-1}}.$$
Multiplying equation (\ref{eq5}) by $\eta_i\ (i=1,2)$ and
integrating by parts, we obtain
$$\int_{\Omega}|\nabla\psi_i|^{p-2}\nabla\psi_i\cdot\nabla\eta_i+
\beta\int_{\partial\Omega}\psi_i^{p-1}\eta_i-\lambda_1\int_{\Omega}\psi_i^{p-1}\eta_i=0,\
\ (i=1,2).$$ It follows
\begin{equation}\label{eq6}
\begin{array}{ll}
&\int_{\Omega}(1+(p-1)(\frac{\psi_2}{\psi_1})^p)|\nabla\psi_1|^p+
\int_{\Omega}(1+(p-1)(\frac{\psi_1}{\psi_2})^p)|\nabla\psi_2|^p\\[3mm]
&-\int_{\Omega}(p(\frac{\psi_2}{\psi_1})^{p-1}|\nabla\psi_1|^{p-2}
+p(\frac{\psi_1}{\psi_2})^{p-1}|\nabla\psi_2|^{p-2})\nabla\psi_1\cdot\nabla\psi_2=0.
\end{array}
\end{equation}
Noticing that $\nabla(\ln\psi_i)=\frac{\nabla\psi_i}{\psi_i}$,
(\ref{eq6}) can be rewritten as
\begin{equation}\label{eq7}
\begin{array}{ll}
&\int_{\Omega}(\psi_1^p+(p-1)\psi_2^p)|\nabla\ln\psi_1|^p+(\psi_2^p+(p-1)\psi_1^p)|\nabla
\ln\psi_2|^p\ \ \ \ \ \\[3mm]
&=p\int_{\Omega}(\psi_2^p|\nabla \ln\psi_1|^{p-2}+\psi_1^p|\nabla
\ln\psi_2|^{p-2})\nabla\ln\psi_1\cdot\nabla\ln\psi_2
\end{array}
\end{equation}
Hence
\begin{equation}\label{eq8}
\begin{array}{ll}
\int_{\Omega}(\psi_1^p-\psi_2^p)(|\nabla \ln\psi_1|^p-|\nabla
\ln\psi_2|^p)\\
= p\int_{\Omega}\psi_2^p|\nabla \ln\psi_1|^{p-2}(\nabla
\ln\psi_1)\cdot(\nabla \ln\psi_2-\nabla \ln\psi_1)\\
-p\int_{\Omega}\psi_1^p|\nabla \ln\psi_2|^{p-2}(\nabla
\ln\psi_2)\cdot(\nabla \ln\psi_2-\nabla \ln\psi_1).
\end{array}
\end{equation}
Observing that (see \cite{sLind})
\begin{equation}\label{eq9}
|\xi_2|^p-|\xi_1|^p\geq p
|\xi_1|^{p-2}\xi_1\cdot(\xi_2-\xi_1)+C(p)\frac{|\xi_2-\xi_1|^p}{2^p-1},\
\ \forall\ \xi_1,\xi_2\in R^n,
\end{equation}
we obtain
\begin{equation}\label{eq10}
\begin{array}{ll}
|\nabla \ln\psi_1|^p-|\nabla \ln\psi_2|^p\\
\geq p |\nabla
\ln\psi_2|^{p-2}(\nabla \ln\psi_2)\cdot(\nabla \ln\psi_1-\nabla
\ln\psi_2)+C_1(p)\frac{|\nabla \ln\psi_1-\nabla
\ln\psi_2|^p}{2^p-1},
\end{array}
\end{equation}
and
\begin{equation}\label{eq11}
\begin{array}{ll}
|\nabla \ln\psi_2|^p-|\nabla \ln\psi_1|^p\\
\geq p |\nabla \ln\psi_1|^{p-2}(\nabla \ln\psi_1)\cdot(\nabla
\ln\psi_2-\nabla \ln\psi_1)+C_2(p)\frac{|\nabla \ln\psi_2-\nabla
\ln\psi_1|^p}{2^p-1},
\end{array}
\end{equation}
where $C(p)$, $C_1(p)$ and $C_2(p)$ are positive constants depend
only on $p$.

From (\ref{eq8}), (\ref{eq10}) and (\ref{eq11}), we deduce
$$-\frac{C_1(p)+C_2(p)}{2^p-1}\int_{\Omega}(\frac{1}{\psi_2^p}+\frac{1}{\psi_1^p})
|\nabla\ln\psi_1-\nabla\ln\psi_2|^p\geq 0.$$ This implies that
$\nabla (\ln\psi_1-\ln\psi_2)=0$,  namely, $\psi_2=C\psi_1$. This
completes the proof of Proposition 2.2.

\vskip 0.2cm

{\em {\bf Corollary 2.3.} If $\Omega=B(0)$ is a ball, then the
first eigenfunction $\psi$ of problem (\ref{eq1}) is radially
symmetry, that is, $\psi(x)=\psi(r)$ with $r=|x|$.}

{\bf Proof.} The conclusion of Corollary 2.3 comes immediately
from the simplicity of $\lambda_1(\Omega)$ and the rotational
invariance of problem (\ref{eq1}).

\vskip 0.2cm

To state our next proposition of the first eigenfunction, we need
the following two lemmas which were proved in the appendix of
\cite{sSaka}.

\vskip 0.2cm

{\em {\bf Lemma 2.4.}(Weak comparison principle)\ Let $\Omega\in
R^N$ be a bounded domain with smooth boundary $\partial\Omega$.
Let $u_1, u_2\in W^{1,p}(\Omega)$ satisfy
$$-{\rm{div}}(|\nabla u_1|^{p-2}\nabla u_1)\leq-{\rm{div}}(|\nabla u_2|^{p-2}\nabla
u_2)\ \ \ \ \mbox{in}\ \ \Omega$$ in weak sense. Then $u_1\leq
u_2$ on $\partial\Omega$ implies $u_1\leq u_2$ in $\Omega$.}

\vskip 0.2cm

{\em {\bf Lemma 2.5.}(Hopf's lemma)\  Let $\Omega\in R^N$ be a
bounded domain with smooth boundary $\partial\Omega$. Let $u\in
C^1(\overline{\Omega})$ satisfy
\begin{equation}
\left\{\begin{array}{ll}\label{eq12}
 -{\rm{div}}(|\nabla
u|^{p-2}\nabla u)\geq
0 & x\in\Omega,\\
u>0 & x\in\Omega.
\end{array}
\right.
\end{equation}
If $u=0$ at $x_0\in\partial\Omega$, then $\frac{\partial
u}{\partial\nu}(x_0)<0$, where $\nu$ denotes the unit outward
vector normal to $\partial\Omega$.}

\vskip 0.2cm

We also need the following strong maximum principle which is a
special case of Theorem 1.1 in \cite{sSerrin}.

\vskip 0.2cm

{\em {\bf Lemma 2.6.}(Strong maximum principle)\ If $u\in
C^1(\Omega)$ satisfies the following inequalities in weak sense
\begin{equation}
\left\{\begin{array}{ll}\label{eq13}
 -{\rm{div}}(|\nabla u|^{p-2}\nabla u)\geq
0 & x\in\Omega,\\
u\geq 0 & x\in\Omega.
\end{array}
\right.
\end{equation}
Then, $u(x_0)=0$ for some $x_0\in\Omega$ implies $u(x)\equiv 0$ in
$\Omega$.}

\vskip 0.2cm

{\em {\bf Proposition 2.7.}\ Let $B_R(0)$ be a ball in $R^N$ with
radius $R$ and center $0$. If $\psi(x)=\psi(r)$ denotes the first
eigenfunction of problem (\ref{eq1}) on $B_R(0)$, then
$\psi'(r)<0$ for any $0<r\leq R$.}

\vskip 0.2cm

{\bf Proof:}\ For any fixed $r_0\in (0,R)$, we have

\[\left\{\begin{array}{ll}
-{\rm{div}}(|\nabla\phi(r_0)|^{p-2}\nabla
\phi(r_0))\leq-{\rm{div}}(|\nabla
\phi(x)|^{p-2}\nabla\phi(x))& x\in B_{r_0}(0)\\
\phi(x)=\phi(r_0) & x\in\partial B_{r_0}(0).
\end{array}
\right.\]
Hence, by lemma 2.4, we have
$$\phi(x)\geq\phi(r_0)\ \ \ \ \ \ x\in B_{r_0}(0).$$
Since $\phi(x)$ is not a constant, it follows from lemma 2.6 that
$$\phi(x)>\phi(r_0)\ \ \ \ \ \ x\in B_{r_0}(0).$$
Let $w(x)=\phi(x)-\phi(r_0)=\phi(r)-\phi(r_0)$. Then, $w(x)$
satisfies
\[\left\{\begin{array}{ll}
-{\rm{div}}(|\nabla w(x)|^{p-2}\nabla w(x))=\lambda_1\phi(x)>0 & x\in B_{r_0}(0)\\
w(x)>0 & x\in B_{r_0}(0)\\
w(x)=0 & x\in\partial B_{r_0}(0).
\end{array}
\right.\]
Consequently, lemma 2.5 implies that $\phi'(r_0)<0$.
Noting that $r_0$ is arbitrary, the conclusion of proposition 2.7
then follows.

\vskip 0.2cm

We conclude this section with the following proposition which is
essential for the proof of Theorem 1.1.

\vskip 0.2cm

{\em {\bf Proposition 2.8.}\ Let $B_R(0)$ be a ball in $R^N$ with
radius $R$ and center $0$. Let $\psi(x)=\psi(r)$ denote the first
eigenfunction of problem (\ref{eq1}) on $B_R(0)$. If
$g(r)=|\psi'(r)|/\psi(r)$, then $g'(r)>0$ for $0<r<R$, and
$g(r)\leq\beta^{\frac{1}{p-1}}$ for any $r\in [0,R]$.}

\vskip 0.2cm

{\bf Proof:}\ It follows from Proposition 2.7 and the standard
regularity theory of elliptic equations that $\psi\in
C^{\infty}(B_R(0)\setminus\{0\})$. Consequently, $0<g\in
C^{\infty}(0,R)$. Now, we compute
\begin{equation}\label{eq14}
g'=(-\frac{\psi'}{\psi})'=-\frac{\psi''}{\psi}+g^2,
\end{equation}
\begin{equation}\label{eq15}
g''=-\frac{\psi'''}{\psi}+3gg'-g^3.
\end{equation}
From the equation satisfied by $\psi$, we have
\begin{equation}\label{eq16}
-(p-1)|\psi'|^{p-2}\psi''-\frac{N-1}{r}|\psi'|^{p-2}\psi'=\lambda_1\psi^{p-1}.
\end{equation}
It follows that
\begin{equation}\label{eq17}
-(p-1)\psi''-\frac{N-1}{r}\psi'=\lambda_1\frac{\psi}{g^{p-2}}.
\end{equation}
Differentiating the above equation, we obtain
\begin{equation}\label{eq18}
-(p-1)\frac{\psi'''}{\psi}-\frac{N-1}{r}\frac{\psi''}{\psi}-\frac{N-1}{r^2}g=-\frac{\lambda_1}{g^{p-3}}-
\frac{(p-2)\lambda_1g'}{g^{p-1}}.
\end{equation}
Since $-\frac{\psi''}{\psi}=g'-g^2$, it follows from (\ref{eq17})
that
\begin{equation}\label{eq19}
\begin{array}{ll}
\lambda_1 &=-(p-1)g^{p-2}\frac{\psi''}{\psi}+\frac{N-1}{r}g^{p-1}\\
          &=(p-1)g^{p-2}g'-(p-1)g^p+\frac{N-1}{r}g^{p-1}.
\end{array}
\end{equation}
Hence
$$-\frac{\lambda_1}{g^{p-3}}=-(p-1)gg'+(p-1)g^3-\frac{N-1}{r}g^2,$$
and

$$-\frac{\psi'''}{\psi}=g^3+\frac{(N-1)g}{(p-1)r^2}-fg',$$
where
$f=g+\frac{N-1}{(p-1)r}+\frac{(p-2)\lambda_1}{(p-1)g^{p-1}}$.
Substituting this equation into  (\ref{eq15}), we infer that for
any $r\in(0,R)$
$$g''(r)+[f(r)-3g(r)]g'(r)=\frac{(N-1)g(r)}{(p-1)r^2}>0.$$
we claim that $g'\neq 0$ in $(0,R)$. For if there exists
$r_0\in(0,R)$ such that $g'(r_0)=0$, then $g''(r_0)>0$. Hence
$r_0$ is a minimum point of $g$. Since $g\geq 0$ and $g(0)=0$, it
follows from the continuity of $g$ that $g(r_0)=0$. This
contradicts with the fact that $g>0$ in $(0,R)$. Consequently,
$g'$ has definite sign in $(0,R)$. This implies immediately that
$g'>0$ in $(0,R)$. For if $g'\leq 0$ in $(0,R)$, then we have
$g(r)\leq g(0)=0$ in $(0,R)$, a contradiction. Finally, the a
priori estimate $g(r)\leq\beta^{\frac{1}{p-1}}$ for any $r\in
[0,R]$ follows from the facts that $g'>0$ and
$g(R)=\beta^{\frac{1}{p-1}}$. This completes the proof of
Proposition 2.8.

\section*{3 Level Sets Formula of $\lambda_1(\Omega)$}

\setcounter{section}{3}

\setcounter{equation}{0}

\noindent For an open set $U\subset\Omega$, we define the interior
and exterior boundary of $U$ respectively by
$$\partial_I U=\partial U\cap\Omega,\ \ \ \partial_E U=\partial U\cap
\partial\Omega.$$
Then $\partial U=\partial_I U\cup\partial_E U$ is a disjoint
union. For any $\varphi\in C(\bar\Omega)$ and $\varphi(x)\geq 0$,
we define a functional $H_{\Omega}(U,\varphi)$ by
\begin{equation}\label{eq20}
H_{\Omega}(U,\varphi)=\frac{1}{|U|}(\int_{\partial_I U}\varphi\
d\sigma +\int_{\partial_E U} \beta\ d\sigma -(p-1)\int_U
\varphi^{\frac{p}{p-1}}\ dx),
\end{equation}
where $\sigma$ is the $(N-1)$-dimensional Hausdorff measure
defined on $\partial U$ and $|U|$ is the Lebesgue measure of $U$.
Since $\varphi$ is continuous on $\partial_I U$ and $U$, all
integrals in $H_{\Omega}(U,\varphi)$ are well defined. In the
following, we reformulate $\lambda_1(\Omega)$ by
$H_{\Omega}(U,\varphi)$. To this end, we always denote by $\psi$
the first eigenfunction of (\ref{eq1}) and sometimes denote $\psi$
by $\psi_{\Omega}$ when we want to emphasize the dependence of
$\psi$ on the domain $\Omega$ . Furthermore, we choose $\psi$ so
that $\psi>0$ and $\|\psi\|_{L^\infty(\Omega)}=1$. By regularity
results of DiBenedetto \cite{sDiBen} and Tolksdorf
\cite{sTolk,sTolk1}, we know that $\psi$ belongs to
$C^{1,\alpha}(\overline{\Omega})$ for some $0<\alpha<1$. Let
$$m=\min_{x\in \bar{\Omega}}\psi(x).$$
Then, by Hopf's boundary point Lemma we have $m>0$. For any
$t\in(m,1)$, we denote by $U_t$ the level set of $\psi$, that is
$$U_t=\{x\in\Omega;\ \psi(x)>t\},$$
then $U_t$ is open and the interior boundary of $U_t$ is the level
surface
$$S_t=\partial_I U_t=\{x\in\Omega;\ \psi(x)=t\}.$$
Hence, $S_t=\emptyset$, if $t\not\in\ (m,1]$. Bearing all
notations $\psi$, $U_t$ and $S_t$ in mind, we prove

\vskip 0.2cm

{\em {\bf Proposition 3.1.} Let $\lambda_1(\Omega)$ be the first
eigenvalue of problem (\ref{eq1}), $\psi$ be the corresponding
eigenfunction , and $H_{\Omega}(U,\varphi)$ be defined as in
(\ref{eq20}). Then
\begin{equation}\label{eq21}
\lambda_1(\Omega)=H_{\Omega}(U_t,\frac{|\nabla\psi|^{p-1}}{\psi^{p-1}})\
\ \ \ \mbox{for}\ \ t\in(m,1).
\end{equation}}

\vskip 0.2cm

As in \cite{sDan}, second partial derivatives of eigenfunction
will be involved in the proof of Proposition 3.1. However, it is
well known that, in general, the best possible regularity results
of Problem (\ref{eq1}) is $C^{1,\alpha}$. Hence, to prove
Proposition 3.1, we consider the following regularized problem
\begin{equation}\label{eq22}
\left\{
\begin{array}{ll}
-{\rm{div}}[(\varepsilon u_{\varepsilon}^2+|\nabla
u_{\varepsilon}|^2)^{\frac{p-2}{2}}\nabla
u_{\varepsilon}]=\lambda^{\varepsilon}
|u_{\varepsilon}|^{p-2}u_{\varepsilon}-\varepsilon(\varepsilon{u_{\varepsilon}}^2+|\nabla
u_{\varepsilon}|^2)^{\frac{p-2}{2}}u_{\varepsilon},\ \
&x\in\Omega,\\[3mm]
(\varepsilon u_{\varepsilon}^2+|\nabla
u_{\varepsilon}|^2)^{\frac{p-2}{2}}\frac{\partial
u_{\varepsilon}}{\partial\nu}+\beta
u_{\varepsilon}^{p-1}=0,&x\in\partial\Omega.
\end{array}
\right.
\end{equation}
where $1<p<+\infty$, $\nu$ is the outward unit normal of
$\partial\Omega$.

Define $\lambda_1^{\varepsilon}(\Omega)$ by
$$\lambda_1^{\varepsilon}(\Omega)=\inf\{\int_{\Omega}
(\varepsilon u_{\varepsilon}^2+|\nabla
u_{\varepsilon}|^2)^{\frac{p}{2}}+\beta\int_{\partial
\Omega}u_{\varepsilon}^p;\ u_{\varepsilon}\in\mathcal{K}\}.$$
Then, we have

\vskip 0.2cm

{\em{\bf Lemma 3.2.} For any $\varepsilon>0$,
$\lambda_1^{\varepsilon}(\Omega)$ is attained by a positive
function $\psi_{\varepsilon}\in\mathcal{K}$. Moreover, up to a
subsequence, we have $\lim\limits_{\varepsilon\rightarrow 0}
\lambda_1^{\varepsilon}(\Omega) =\lambda_1(\Omega)$, and
$$\lim\limits_{\varepsilon\rightarrow 0}\psi_{\varepsilon}=\widehat{\psi}
\ \ \ \ \mbox{in}\ \ C^1(\overline{\Omega})$$ where
$\lambda_1(\Omega)$ is the first eigenvalue of problem (\ref{eq1})
and $\widehat{\psi}$ is a corresponding eigenfunction with
$||\widehat{\psi}||_{L^p(\Omega)}=1$. }

\vskip 0.2cm

{\bf Proof.} The conclusion that $\lambda_1^{\varepsilon}(\Omega)$
is attained by a positive function $\psi_{\varepsilon}$ can be
proved in the same way as that of Proposition 2.1. To prove the
second part of Lemma 3.2, We first note that $\psi_{\varepsilon}$
is bounded in $W^{1,p}(\Omega)$ when $\varepsilon$ is small
enough. Hence, up to a subsequence, we may assume that
$$\psi_{\varepsilon}\rightharpoonup\tilde{\psi}\ \ \rm{weakly\ in\ }
W^{1,p}(\Omega),\ \ \rm{as}\ \varepsilon\rightarrow 0.$$ Since
$W^{1,p}(\Omega)\hookrightarrow L^p(\Omega)$ is compact, we also
have $\tilde{\psi}\in\mathcal{K}$. Because convex functional is
weakly lower semi-continuous, we have
\begin{equation}\label{eq23}
\int_{\Omega}|\nabla\tilde{\psi}|^p+\beta\int_{\partial\Omega}|\tilde{\psi}|^p
\leq\underline{\lim}_{\varepsilon\rightarrow
0}(\int_{\Omega}|\nabla\psi_{\varepsilon}|^p+\beta\int_{\partial\Omega}|\psi_{\varepsilon}|^p)
\end{equation}
On the other hand, if $\widehat{\psi}$ is the first eigenfunction
of problem (\ref{eq1}) with $||\widehat{\psi}||_{L^p(\Omega)}=1$,
then by the definition of $\lambda_1(\Omega)$ and
$\lambda_1^\varepsilon(\Omega)$, we have
\begin{equation}\label{eq24}
\begin{array}{ll}
\lambda_1(\Omega)&=\int_{\Omega}|\nabla\widehat{\psi}|^p+\beta\int_{\partial\Omega}|\widehat{\psi}|^p
\leq\int_{\Omega}|\nabla\psi_{\varepsilon}|^p+\beta\int_{\partial\Omega}|\psi_{\varepsilon}|^p\\
&\leq\int_{\Omega}(\varepsilon\psi_{\varepsilon}^2+|\nabla\psi_{\varepsilon}|^2)^{\frac{p}{2}}
+\beta\int_{\partial\Omega}\psi_{\varepsilon}^p\\[3mm]
&=\lambda_{\varepsilon}(\Omega)\leq
\int_{\Omega}(\varepsilon\widehat{\psi}^2+|\nabla\widehat{\psi|}^2)^{\frac{p}{2}}+\beta\int_{\partial\Omega}
\widehat{\psi}^p
\end{array}
\end{equation}
Let $\varepsilon\rightarrow 0$ on the both side of the above
inequality, we obtain
\begin{equation}\label{eq25}
\lim_{\varepsilon\rightarrow
0}\lambda_{\varepsilon}(\Omega)=\lim_{\varepsilon\rightarrow
0}\int_{\Omega}|\nabla\psi_{\varepsilon}|^p+\beta\int_{\partial\Omega}|\psi_{\varepsilon}|^p
=\int_{\Omega}|\nabla\widehat{\psi}|^p+\beta\int_{\partial\Omega}|\widehat{\psi}|^p=\lambda_1(\Omega)
\end{equation}
From (\ref{eq23}) and (\ref{eq25}), we infer that
$$\int_{\Omega}|\nabla\tilde{\psi}|^p+\beta\int_{\partial\Omega}|\tilde{\psi}|^p
\leq\int_{\Omega}|\nabla\widehat{\psi}|^p+\beta\int_{\partial\Omega}|\widehat{\psi}|^p=\lambda_1(\Omega).$$
Hence, $\tilde{\psi}$ is a minimizer of $\lambda_1(\Omega)$. This
implies that $\tilde{\psi}=\widehat{\psi}$ due to the simplicity
of $\lambda_1(\Omega)$ and
$||\tilde{\psi}||_{L^p(\Omega)}=||\widehat{\psi}||_{L^p(\Omega)}=1$.
Consequently, $\psi_{\varepsilon}\rightharpoonup\widehat{\psi}$
weakly in $W^{1,p}(\Omega)$ as $\varepsilon\rightarrow 0$.
Finally, by the regularity theory of Tolksdorf \cite{sTolk,sTolk1}
and DiBenedetto \cite{sDiBen}, we know that for
$\varepsilon\in(0,1)$,there exists $\alpha\in(0,1)$ and a positive
constant $C$ independent of $\varepsilon$ such that
$||\psi_{\varepsilon}||_{C^{1,\alpha}}(\overline{\Omega})\leq C$.
Hence, up to a subsequence, $\psi_\varepsilon$ converges to
$\widehat{\psi}$ in $C^1(\overline{\Omega})$. This completes the
proof of Lemma 3.2.

\vskip 0.2cm

{\bf Proof of Proposition 3.1.} For any fixed $t\in(m,1)$, let
$\nu$ denote the outward unit vector normal to $\partial U_t$. If
we denote by $\psi_\varepsilon$ the solution of Problem
(\ref{eq22}) obtained in Lemma 3.2, then by the standard
regularity theory of elliptic equations we know that
$\psi_\varepsilon\in C^{\infty}(\Omega)$. Hence, by divergence
Theorem, we have
\begin{equation}\label{eq26}
\begin{array}{ll}
-\int_{\partial U_t}
\frac{(\varepsilon\psi_{\varepsilon}^2+|\nabla\psi_{\varepsilon}|^2)^{\frac{p-2}{2}}}{\psi_{\varepsilon}^{p-1}}
\frac{\partial\psi_{\varepsilon}}{\partial\nu}d\sigma\\
=-\int_{U_t}{\rm{div}}
(\frac{(\varepsilon\psi_{\varepsilon}^2+|\nabla\psi_{\varepsilon}|^2)^{\frac{p-2}{2}}\nabla\psi_{\varepsilon}}{\psi_{\varepsilon}^{p-1}})dx\\
=-\int_{U_t}\frac{{\rm{div}}((\varepsilon\psi_{\varepsilon}^2+|\nabla\psi_{\varepsilon}|^2)^{\frac{p-2}{2}}\nabla\psi_{\varepsilon})}{\psi_{\varepsilon}^{p-1}}
dx
+(p-1)\int_{U_t}\frac{(\varepsilon\psi_{\varepsilon}^2+|\nabla\psi_{\varepsilon}|^2)^{\frac{p-2}{2}}|\nabla\psi_{\varepsilon}|^2}{\psi_{\varepsilon}^p}dx\\
=\lambda_1^{\varepsilon}(\Omega)|U_t|-\varepsilon\int_{U_t}\frac{(\varepsilon\psi_\varepsilon^2+|\nabla\psi_\varepsilon|^2)^{\frac{p-2}{2}}}{\psi_\varepsilon^{p-2}}
+(p-1)\int_{U_t}\frac{(\varepsilon\psi_{\varepsilon}^2+|\nabla\psi_{\varepsilon}|^2)^{\frac{p-2}{2}}|\nabla\psi_{\varepsilon}|^2}{\psi_{\varepsilon}^p}dx
\end{array}
\end{equation}
Passing to the limit in (\ref{eq26}) as $\varepsilon\rightarrow
0$, we obtain
\begin{equation}\label{eq27}
-\int_{\partial U_t}
\frac{|\nabla\widehat{\psi}|^{p-2}}{\widehat{\psi}^{p-1}}
\frac{\partial\widehat{\psi}}{\partial\nu}d\sigma
=\lambda_1(\Omega)|U_t|
+(p-1)\int_{U_t}\frac{|\nabla\widehat{\psi}|^{p}}{\widehat{\psi}^p}dx
\end{equation}
Since $\lambda_1(\Omega)$ is simple, we have
\begin{equation}\label{eq50}
-\int_{\partial U_t} \frac{|\nabla\psi|^{p-2}}{\psi^{p-1}}
\frac{\partial\psi}{\partial\nu}d\sigma =\lambda_1(\Omega)|U_t|
+(p-1)\int_{U_t}\frac{|\nabla\psi|^{p}}{\psi^p}dx.
\end{equation}
By the boundary condition, we have
$$\beta=-\frac{|\nabla\psi|^{p-2}}{\psi^{p-1}}\frac{\partial\psi}
{\partial \nu},\ \ \ x\in \partial_E U_t,$$ due to $\partial_E
U_t\subset
\partial \Omega$.
Noticing further that $|\nabla\psi|=-\frac{\partial
\psi}{\partial\nu}$ on $S_t$, we obtain from the definitions of
$S_t$ and $\partial_EU_t$ that
\begin{equation}\label{eq28}
-\int_{\partial U_t}\frac{\ \ |\nabla\psi|^{p-2}}{\psi^{p-1}}
\frac{\partial\psi}{\partial
\nu}d\sigma=\int_{S_t}\frac{|\nabla\psi|^{p-2}}{\psi^{p-1}}|\nabla\psi|
d\sigma+\int_{\partial_E U_t}\beta \ d\sigma .
\end{equation}
Now, the conclusion of Proposition 3.1 follows from (\ref{eq50})
and (\ref{eq28}).

\section*{4 The lower bound of $\lambda_1(\Omega)$ }
\setcounter{section}{4}

\setcounter{equation}{0}

\noindent In this section, we give a lower bound of
$\lambda_1(\Omega)$. Let
\begin{equation}\label{eq29}
\mathcal{M}_{\beta}=\{\varphi(x)\in C(\Omega);\varphi(x)\geq0,
\overline{\lim_{x\to z}}\varphi(x)\leq \beta,\ \ \forall
z\in\partial\Omega\}.
\end{equation}
Keep in use the same notations $\psi$, $U_t$ and $S_t$ as in the
previous section. Since $\psi(x)\in C^1(\overline{\Omega})$, it is
easy to see that
$(\frac{|\nabla\psi|}{\psi})^{p-1}\in\mathcal{M}_{\beta}$ if and
only if $\psi$ is a constant on $\partial \Omega$. In fact, if
$\psi$ is a constant on $\partial \Omega$ then $\frac{\partial
\psi}{\partial \nu}=-|\nabla\psi|$ on $\partial \Omega$. Hence
$$\frac{|\nabla\psi|^{p-1}}{\psi^{p-1}}=-\frac{|\nabla\psi|^{p-2}}{\psi^{p-1}}\frac{\partial
\psi}{\partial \nu}=\beta\ {\rm {on}}\ \partial\Omega.\ \ \ \ $$
This implies that
$(\frac{|\nabla\psi|}{\psi})^{p-1}\in\mathcal{M}_{\beta}$.

On the other hand, if
$(\frac{|\nabla\psi|}{\psi})^{p-1}\in\mathcal{M}_{\beta}$, then
$$\frac{|\nabla\psi|^{p-1}}{\psi^{p-1}}\leq\beta=-\frac{|\nabla\psi|^{p-2}}{\psi^{p-1}}
\frac{\partial\psi}{\partial
\nu}\leq\frac{|\nabla\psi|^{p-1}}{\psi^{p-1}}\ \ \forall\
x\in\partial\Omega.$$ Hence
$\frac{\partial\psi}{\partial\nu}=-|\nabla\psi|$ for all
$x\in\partial\Omega$, which implies that $\psi$ is a constant on
$\partial\Omega$.

The main results of this section can be stated as

\vskip 0.2cm

{\em{\bf Theorem 4.1.} For every $\varphi\in \mathcal{M}_{\beta}$,
there exists a set $I\subset [0,1]$ with positive measure such
that
\begin{equation}\label{eq30}
\lambda_1(\Omega)\geq H_{\Omega}(U_t,\varphi)\ \ \ \ \mbox{for
all}\ \ \ t\in I .
\end{equation}}

\vskip 0.2 cm

{\em{\bf Theorem 4.2.} Let $\varphi(x)\in\mathcal{M}_{\beta}$, and
$\psi$ be the first eigenfunction of problem (\ref{eq1}). If
$\varphi\neq\frac{|\nabla\psi|^{p-1}}{\psi^{p-1}}$, then there
exists a set $I\subset[m,1]$ with positive measure such that
$$H_{\Omega}(U_t,\varphi)<\lambda_1(\Omega) \ \ \mbox{for all}\ \ t\in I.$$}

\vskip 0.2cm

To prove Theorems, we prove some lemmas first. For any given
$\varphi$ and $\varphi\geq0$, let
\begin{equation}\label{eq31}
\omega(x):=\varphi(x)-\frac{|\nabla\psi|^{p-1}}{\psi^{p-1}},\ \
x\in\Omega.
\end{equation}
Then we have

\vskip 0.2cm

{\em{\bf Lemma 4.3.} For any $\varphi\in \mathcal{M}_{\beta}$, let
$\omega$ be defined as (\ref{eq31}). Then for any $\varepsilon>0$
there exists $\delta>0$ such that $\omega(x)\leq\varepsilon$ for
all $x\in\Omega$ with ${\rm{dist}}(x,\partial\Omega)<\delta$.}

\vskip 0.2cm

{\bf Proof.} Since $\frac{|\nabla\psi|^{p-1}}{\psi^{p-1}}$ is
continuous on the compact set $\overline{\Omega}$, we have that
for any fixed $\varepsilon>0$, there exists $\delta_0>0$ such that
\begin{equation}\label{eq32}
|\frac{|\nabla\psi(x)|^{p-1}}{\psi(x)^{p-1}}-\frac{|\nabla\psi(z)|^{p-1}}
{\psi(z)^{p-1}}|<\frac{\varepsilon}{2}\ \ \  {\rm{for\ any}}\ \ x,\
z\in\overline{\Omega},\ |x-z|<\delta_0.
\end{equation}
For any fixed $z\in\partial\Omega$ and the above fixed
$\varepsilon$, by the assumption that $\overline{\lim}_{x\to
z}\varphi\leq\beta$, we can choose $r_z>0$ such that
$$\sup_{x\in B(z,r_z)\cap\Omega}\varphi(x)\leq\beta+\frac{\varepsilon}{2},$$
that is
\begin{equation}\label{eq33}
\varphi(x)-\beta\leq\frac{\varepsilon}{2},\ \ {\rm{for\ all}}\ x\in
B(z,r_z)\cap\Omega,
\end{equation}
where $B(z,r_z)$ denotes the ball with radius $r_z$ and center $z$
. Since the set $\{B(z,r_z),\ z\in\partial\Omega\}$ of balls form
an open cover of the compact set $\partial\Omega$, we can select a
finite sub-cover $\{B(z_i,r_i)\}_{i=1}^n$ with $r_i=r_{z_i}$. Let
$\delta\leq \min\{r_1,r_2,\cdots,r_n,\delta_0\}$ be so small that
$x\in\bigcup\limits_{i=1}^nB(z_i,r_i)$ whenever $x\in\Omega$
satisfying dist$(x,\partial\Omega)<\delta$. Then, for any
$x\in\Omega$ with dist$(x,\partial\Omega)<\delta$, there exists
$i_0\in\{1,2,\cdots,n\}$ such that $x\in B(z_{i_0},r_{i_0})$. By
the boundary condition, we have
\begin{equation}\label{eq60}
\beta=-\frac{|\nabla\psi|^{p-2}}{\psi^{p-1}}\frac{\partial
\psi}{\partial
\nu}(z_{i_0})\leq\frac{|\nabla\psi|^{p-1}}{\psi^{p-1}}(z_{i_0}).
\end{equation}
It follows from (\ref{eq32}), (\ref{eq33}) and (\ref{eq60}) that
for any $x\in\Omega$ with dist$(x,\partial\Omega<\delta$, we have
\begin{equation}\label{eq34}
\begin{array}{ll}
\omega(x)=\varphi(x)-\frac{|\nabla\psi(x)|^{p-1}}{\psi(x)^{p-1}}
&=\varphi(x)-\beta+\beta-\frac{|\nabla\psi(x)|^{p-1}}{\psi(x)^{p-1}}\\
&\leq\frac{\varepsilon}{2}+\frac{|\nabla\psi(z_{i_0})|^{p-1}}{\psi(z_{i_0})^{p-1}}
-\frac{|\nabla\psi(x)|^{p-1}}{\psi(x)^{p-1}}\\
&\leq\frac{\varepsilon}{2}+\frac{\varepsilon}{2}=\varepsilon.
\end{array}
\end{equation}
This is just the desired conclusion of Lemma 4.3.

\vskip 0.2cm

{\em{\bf Lemma 4.4.} Suppose that $\varphi\in C(\Omega)$ is
non-negative such that $\varphi\in L^1(U)$ for every open set
$U\subset\Omega$. Let $\omega$ be defined as (\ref{eq31}). Set
$$F(t):=\int_t^1\frac{1}{\tau}\int_{S_{\tau}}\omega d\sigma d\tau,\ \ {\rm{for}}\ t\in(m,1).$$
Then $F$ is absolutely continuous on $(\varepsilon,1)$ for all
$\varepsilon\in(0,1)$ and
$$\frac{d}{dt}F(t)=-\frac 1t\int_{S_t}\omega d\sigma,$$
for almost all $t\in(0,1)$. }

\vskip 0.2cm

{\bf Proof.} Fix $\varepsilon\in(0,1)$. By the assumption
$\varphi\in C(\Omega)\cap L^1(U_{\varepsilon})$ and the co-area
formula, we have
$$\int_{\varepsilon}^1\frac{1}{\tau}\int_{S_{\tau}}\varphi\ d\sigma
d\tau=\int_{U_\varepsilon}\frac{\varphi}{\psi}|\nabla\psi|\
dx<\infty$$ and
$$\int_{\varepsilon}^1\frac{1}{\tau}\int_{S_{\tau}}\frac{|\nabla\psi|^{p-1}}{\psi^{p-1}}\ d\sigma\
d\tau=\int_{U_\varepsilon}\frac{|\nabla\psi|^p}{\psi^p}\
dx<\infty.$$ Let
$$f(\tau):=\frac{1}{\tau}\int_{S_{\tau}}\omega\ d\sigma,$$
Then $f(\tau)\in L^1((\varepsilon,1))$, thus
$F(t)=\int_t^1f(\tau)\ d\tau$ is absolutely continuous on
$(\varepsilon, 1)$ and differentiable almost everywhere. Moreover
$$F'(t)=-f(t)=-\frac{1}{t}\int_{S_{t}}\omega\ d\sigma.$$
This completes the proof of Lemma 4.4.

\vskip 0.2cm

To state our next Lemma, we recall more regularity results of the
first eigenfunction $\psi$. By the boundary condition and the
Hopf's boundary point Lemma, we know that $\psi(x)>0$ for any
$x\in\partial\Omega$. Consequently $|\nabla\psi|(x)>0$ for any
$x\in\partial\Omega$. Since $\partial\Omega$ is compact and
$|\nabla\psi|\in C(\overline{\Omega})$, it is easy to prove that
there exists positive number $\alpha$ and a neighborhood $N$ of
$\partial\Omega$ in $\Omega$ such that
$|\nabla\psi|(x)\geq\alpha>0$ for any $x\in N$. This implies that
p-Laplacian is uniformly elliptic in $N$. Hence, by the interior
regularity theorem of elliptic equations, we know that $\psi\in
C^{\infty}(N)$. If we let
$m=\min\{\psi(x);x\in\overline{\Omega}\}$ and
$K=\{x\in\Omega;\psi(x)=m\}$, then by strong maximum principle we
know that $K\subset\partial\Omega$. Noticing furthermore that $K$
is compact, there exists $t_0\in(m,1)$ small enough such that
$$S_t\subset N\ \ \ \ \mbox{for any}\ \ \ t\leq t_0.$$
An argument similar to that used by Daners in \cite{sDan} implies
the following lemma since all computations in \cite{sDan} are
local.

\vskip 0.2cm

{\em{\bf Lemma 4.5.} Let $U_t$ and $S_t$ be defined as in section
3. Then $U_t$ is a Lipschitz domain, moreover, there exist
$t_1\in(m,t_0)$ and a constant $C>0$ independent of $t$ such that
$\sigma(S_t)\leq C\sigma(\partial \Omega)$ for all $t\in
(m,t_1)$.}

\vskip 0.2cm

{\bf Proof of Theorem 4.1.} We give a proof by contradiction.
Suppose that there exists $\varphi\in\mathcal{M}_{\beta}$ such that
\begin{equation}\label{eq35}
\lambda_1(\Omega)<H_{\Omega}(U_t,\varphi)\ \ {\rm{for\ almost\
all}}\ t\in(m,1).
\end{equation}
Let $\omega$ be defined as (\ref{eq31}) and $F(t)$ be defined as
in Lemma 4.4, that is
\begin{equation}\label{eq65}
\omega(x):=\varphi(x)-\frac{|\nabla\psi|^{p-1}}{\psi^{p-1}},\ \
x\in\Omega.
\end{equation}
and
$$F(t):=\int_t^1\frac{1}{\tau}\int_{S_{\tau}}\omega d\sigma d\tau,\ \ {\rm{for\ all}}\ t\in(m,1).$$
 Then by (\ref{eq35}), the definition of
$H_{\Omega}(U_t,\varphi)$ and Proposition 3.1, we have
\begin{equation}\label{eq36}
\int_{S_t}\omega\ d\sigma
-(p-1)\int_{U_t}(\varphi^{\frac{p}{p-1}}-\frac{|\nabla\psi|^p}{\psi^p})\
dx=|U_t|[H_{\Omega}(U_t,\varphi)-\lambda_1(\Omega)]>0.
\end{equation}
By Taylor's expansion, there holds
\begin{equation}\label{eq37}
\begin{array}{ll}
\varphi^{\frac{p}{p-1}}-\frac{|\nabla\psi|^p}{\psi^p}=
(\frac{|\nabla\psi|^{p-1}}{\psi^{p-1}}+\omega)^{\frac{p}{p-1}}
-\frac{|\nabla\psi|^p}{\psi^p}\\
=\frac{p}{p-1}(\frac{|\nabla\psi|^{p-1}}{\psi^{p-1}})^{\frac{1}{p-1}}\omega
+\frac 12\frac 1{p-1}\xi^{\frac p{p-1}-2}\omega^2,
\end{array}
\end{equation}
where $\xi$ is a nonnegative function with value between $\varphi$
and $\frac{|\nabla\psi|^{p-1}}{\psi^{p-1}}$.

From (\ref{eq36}), (\ref{eq37}), the co-area formula and the
definition of $F(t)$, we obtain
$$
\int_{S_t}\omega\
d\sigma>p\int_{U_t}\frac{|\nabla\psi|}{\psi}\omega\ dx=
p\int_t^1\int_{S_\tau}\frac 1{\tau}\omega\ d\sigma d\tau=pF(t)$$
for almost all $t\in(m,1)$.

It follows from Lemma 4.4 and the above inequality that
$$\frac d{dt}(t^pF(t))=-t^pf(t)+pt^{p-1}F(t)=t^{p-1}(-\int_{S_t}\omega\
d\sigma+pF(t))<0$$ for almost all $t\in(m,1)$.

Hence, the function $t^pF(t)$ is strictly decreasing on $(m,1)$.
Since $F(1)=0$ and $F(t)$ is continuous on $(m,1)$, there exists
$\eta>0$ and $t_2\in(m,1)$ such that $F(t)>\eta$ for
$t\in(m,t_2].$ On the other hand, by Lemma 4.5, there exists
$t_3\in (m,t_2]$ and a constant $C>0$ such that $\sigma(S_t)\leq
C\sigma(\partial\Omega)$ for $t\in (m,t_3)$. Set
$$\varepsilon_0=\frac{\eta}{C\sigma(\partial\Omega)}.$$
For this fixed $\varepsilon_0$, it follows from Lemma 4.3 that
there exists $\delta_0>0$ such that $\omega(x)\leq\varepsilon_0$
for any $x\in\Omega$ with dist$(x,\partial\Omega)<\delta_0$.
Noticing that $\psi$ attains its strict minimum on
$\partial\Omega$, we can choose $0<t_4<t_3$ so small that
dist$(x,\partial\Omega)<\delta_0$ for any $x\in S_t$ and
$t\in(m,t_4)$. Hence, for any $t\in(m,t_4)$, there holds
$$p\eta<pF(t)<\int_{S_t}\omega\ d\sigma\leq\varepsilon_0\sigma(S_t)
\leq\varepsilon_0 C\sigma(\partial\Omega)\leq\eta$$ which is a
contradiction. Thus we complete the proof of Theorem 4.1.

\vskip 0.2cm

{\bf Proof of Theorem 4.2.} We give a proof by contradiction.
Assume that $\varphi\neq\frac{|\nabla \psi|^{p-1}}{\psi^{p-1}}$
and that
\begin{equation}\label{eq38}
H_{\Omega}(U_t,\varphi)\geq\lambda_1(\Omega),\ \ {\rm{for\ almost\
all }}\ t\in(m,1).
\end{equation}
Similar to the proof of Theorem 4.1, by the definition of
$H_{\Omega}(U_t,\varphi)$ and Proposition 3.1, we have
\begin{equation}\label{eq39}
\int_{S_t}\omega\ d\sigma\geq pF(t)+\frac
12\frac{p}{p-1}\int_{U_t}\xi^{\frac{2-p}{p-1}}\omega^2dx\ \ \
\mbox{for almost all}\ \ t\in(m,1),
\end{equation}
and
$$\frac d{dt}(t^pF(t))\leq-\frac
12\frac{p}{p-1}\int_{U_t}\xi^{\frac{2-p}{p-1}}\omega^2dx\leq 0,
 \ \ \mbox{for almost all}\ \ t\in(m,1),$$
where $\xi$ is a nonnegative function with value between $\varphi$
and $\frac{|\nabla \psi|^{p-1}}{\psi^{p-1}}$. Hence $t^pF(t)$ is
nonincreasing in $(m,1)$. Since $\omega(x)\in C(\Omega)$,
$\omega(x)\not\equiv 0$ and
$\bigcup\limits_{t\in(m,1)}U_t=\Omega$, there exists $t_0\in(m,1)$
such that
\begin{equation}\label{eq40}
\int_{U_{t_0}}\xi^{\frac{2-p}{p-1}}\omega^2dx> 0
\end{equation}
Moreover, if $t_1,t_2\in(m,1)$ satisfy $t_1<t_2$, then we have
$U_{t_2}\subset U_{t_1}$. Hence, the map
$$t\mapsto\int_{U_t}\xi^{\frac{2-p}{p-1}}\omega^2dx$$
is non-increasing in $(m,1)$ and
$\int_{U_1}\xi^{\frac{2-p}{p-1}}\omega^2dx=0$ due to
$U_1=\emptyset$.

Let
$$t^*=\sup\{t\in(m,1),\ \ \int_{U_t}\xi^{\frac{2-p}{p-1}}\omega^2dx>0\}.$$
From (\ref{eq40}), we know that $t^*\in(m,1]$ and thus $t^pF(t)$
is strictly decreasing on $(m,t^*)$ and non-increasing on
$[t^*,1]$, similar to the proof of Theorem 4.1, there exists
$t_3\in(m,t^*)$ such that for any $t\in(m, t_3)$,
$$p\eta<pF(t)<\int_{S_t}\omega\ d\sigma\leq\varepsilon\sigma(S_t)\leq\eta,$$
which is a contradiction. Hence, we complete the proof of Theorem
4.2.

\section*{5 Proof of Theorem 1.1}

\setcounter{section}{5}

\setcounter{equation}{0}

\noindent This section devotes to prove Theorem 1.1. To this end,
we denote by $\lambda_1(\Omega)$ the first eigenvalue of problem
(\ref{eq1}) on the domain $\Omega$ and $\psi_{\Omega}$ denotes its
corresponding eigenfunction. Furthermore, $B=B_R(0)$ denotes the
ball with radius $R$ and center $0$ such that $|B|=|\Omega|$. Let
$U_t$ be the level set and $S_t$ be the level surface of
$\psi_{\Omega}$ at level $t$ defined in section 2, and
$B_{r(t)}(0)$ be the ball with radius $r(t)$ and center $0$ such
that $|B_{r(t)}(0)|=|U_t|$. Define
$$\Phi_B(x)=\frac{|\nabla\psi_{B}(x)|^{p-1}}{\psi_{B}^{p-1}(x)}\ \ \ \mbox{for}\ \ \ x\in B_R(0).$$
By Corollary 2.3, $\Phi_B$ is radially symmetry. So, we only need
to consider the radial function
$$G(r)=\Phi_B(|x|)=\frac{|\psi_B'(r)|^{p-1}}{\psi_B^{p-1}(r)}=g^{p-1}(r)\ \ \ \mbox{for}\ \ \ r\in(0,R)$$
where $g(r)$ is the function defined in Proposition 2.8. Then by
Proposition 2.8, we know that $G(r)$ is strictly increasing in
$(0,R)$. Consequently, $G(r)\leq G(R)=\beta$ for any $r\in [0,R]$.
we construct our test function as the following.

For any $t\in(m,1)$ and $x\in S_t$, we set
$$\Phi(x)=G(r(t)).$$
It is obvious that $\Phi$ is well defined since $\Omega$ is a
disjoint union of $S_t$, $t\in(m,1]$. Moreover,
$\Phi\in\mathcal{M}_\beta(\Omega)$ due to $\Phi$ is continuous and
$\Phi(x)\leq\beta$ for all $x\in\overline{\Omega}$. It is also not
too difficult to see that
\begin{equation}\label{eq41}
\int_{U_t}|\Phi|^{\frac{p}{p-1}}dx=\int_{B_{r(t)}}\Phi_B^{\frac{p}{p-1}}dx.
\end{equation}
Since by the construction the level sets of $\Phi$ and $\Phi_B$
have the same measure. Now, we are in a position to prove Theorem
1.1.

\vskip 0.2cm

{\bf Proof of Theorem 1.1.}\ Since $\Phi\in\mathcal{M}_{\beta}$,
we conclude from Theorem 4.1 that there exist a set
$I\subset(m,1)$ with positive measure such that
\begin{equation}\label{eq42}
\lambda_1(\Omega)\geq H_{\Omega}(U_{t},\Phi)\ \ \ \mbox{for\ all}\
\ t\in I .
\end{equation}
Noticing that $\sigma(\partial B_{r_t})\leq\sigma(\partial U_t)$
for all $t\in(m,1]$, and $\Phi(x)=G(r(t))\leq\beta$ when $x\in
S_t$, we have
\begin{equation}\label{eq43}
\begin{array}{ll}
\int_{\partial B_{r(t)}}\Phi_B(x)d\sigma=G(r(t))\sigma(\partial
B_{r(t)})&\leq G(r(t))\sigma(\partial U_t)\\
&=G(r(t))(\int_{S_t}d\sigma +\int_{\partial_EU_t}\ d\sigma)\\
&\leq\int_{S_t}\Phi d\sigma +\int_{\partial_E U_t}\beta d\sigma.
\end{array}
\end{equation}
Hence, from (\ref{eq41}), (\ref{eq43}) and the definitions of
$H_{B}(B_{r(t)}(0),\Phi_B)$ and $H_{\Omega}(U_t,\Phi)$, we have
\begin{equation}\label{eq44}
H_{B}(B_{r(t)}(0),\Phi_B)\leq H_{\Omega}(U_t,\Phi)\ \ \ \ \forall
t\in(m,1).
\end{equation}
Since, by Proposition 3.1, we have
$\lambda_1(B)=H_{B}(B_{r(t)}(0),\Phi_B)$ for any $t\in(m,1)$, it
follows from (\ref{eq42}) and (\ref{eq44}) that
$$\lambda_1(\Omega)\geq\lambda_1(B).$$
This completes the proof of Theorem 1.1.

\section*{6 Proof of Theorem 1.2}

\setcounter{section}{6}

\setcounter{equation}{0}

This section devotes to prove Theorem 1.2. To this end, we keep in
use of all notations in section 5, and prove some lemmas first.

\vskip 0.2cm

{\em{\bf Lemma 6.1.} Suppose that $\Omega$ satisfies that
$\lambda_1(\Omega)=\lambda_1(B)$ with $|\Omega|=|B|$. Then
$$\Phi=\frac{|\nabla
\psi_{\Omega}|^{p-1}}{\psi_{\Omega}^{p-1}}\ \ \ \ \mbox{and}\ \ \
H_{\Omega}(U_t,\Phi)=\lambda_1(B),$$ for almost all $ t\in(m,1)$.}

\vskip 0.2cm

{\bf Proof.} If $\lambda_1(\Omega)=\lambda_1(B)$, then by
Proposition 3.1 and (5.4), we have
$\lambda_1(\Omega)=\lambda_1(B)=H_B(B(r(t)),G)\leq
H_{\Omega}(U_t,\Phi)$ for almost all $t\in(m,1)$. Hence by Theorem
4.2, $\Phi=\frac{|\nabla
\psi_{\Omega}|^{p-1}}{\psi_{\Omega}^{p-1}}$. Again, by Proposition
3.1, we obtain
$H_{\Omega}(U_t,\Phi)=\lambda_1(\Omega)=\lambda_1(B)$, for almost
all $t\in(m,1)$.

\vskip 0.2cm

{\em{\bf Lemma 6.2.} Let $\psi_{\Omega}$ be the eigenfunction
corresponding to the first eigenvalue $\lambda_1(\Omega)$, and
$U_t$ be the level set of $\psi_{\Omega}$. Then
$H_{\Omega}(U_t,\Phi)=\lambda_1(B)$ if and only if $U_t$ is a ball
and $\sigma(\partial_EU_t)=0$.}

\vskip 0.2cm

{\bf Proof.} It follows from Proposition 3.1 that
$\lambda_1(B)=H_B(B(r(t)),G)$ for all $t\in(m,1)$. By the
construction of $G$ and $\Phi$, we know that the level sets of $G$
and $\Phi$ have the same measure. Hence
$$\int_{U_t}|\Phi|^{\frac p{p-1}} dx=\int_{B_{r(t)}}|G|^{\frac p{p-1}}
dx,\ \ {\rm{for\ all}}\ t\in(m,1).$$ Using the definitions of
$H_{\Omega}(U,\varphi)$ and $\Phi$, we have
\begin{equation}\label{eq45}
\begin{array}{ll}
H_{\Omega}(U_t,\Phi)&=\frac{1}{|U_t|}(\int_{\partial_I U_t}\Phi
d\sigma +\int_{\partial_E U_t} \beta d\sigma -(p-1)\int_{U_t}
\Phi^{\frac{p}{p-1}} dx),\\[3mm]
&=\frac{1}{|B_{r(t)}|}[G(r(t))\sigma(S_t)+\beta
\sigma(\partial_EU_t)-(p-1)\int_{B_{r(t)}} G^{\frac{p}{p-1}}dx].
\end{array}
\end{equation}
If $U_t$ is a ball and $\sigma(\partial_EU_t)=0$, then
$\sigma(S_t)=\sigma(\partial B_{r(t)})$ and
\begin{equation}\label{eq46}
\begin{array}{ll}
H_{\Omega}(U_t,\Phi)=\frac{1}{|B_{r(t)}|}[G(r(t))\sigma(\partial
B_{r(t)}) -(p-1)\int_{B_{r(t)}}
G^{\frac{p}{p-1}}dx]\\
=H_B(B_{r(t)},G)=\lambda_1(B).
\end{array}
\end{equation}
Conversely, if $H_{\Omega}(U_t,\Phi)=\lambda_1(B)$, then for this
$t$,

$$G(r(t))\sigma(S_t)+\beta
\sigma(\partial_EU_t)=G(r(t))\sigma(\partial B_{r(t)}).$$ Noticing
that $S_t=\partial_I U_t=\partial U_t-\partial_E U_t$, we have
$$\sigma(\partial_EU_t)(\beta-G(r(t)))=G(r(t))(\sigma(\partial B_{r(t)})-\sigma(\partial
U_t)).$$ This is only possible when $\sigma(\partial_EU_t)=0$ and
$\sigma(\partial B_{r(t)})=\sigma(\partial U_t)$, since
$0<G(r(t)))<\beta$ for all $t\in(m,1)$ and $|B_{r(t)}|=|U_t|$
implies $\sigma(\partial B_{r(t)})\leq\sigma(\partial U_t)$. But
we know that the ball is the unique minimizer of the isoperimetric
inequality. Hence, $U_t=B_{r(t)}+z$ for some $z\in R^N$. This
completes the proof of Lemma 6.2.

\vskip 0.2cm

{\em{\bf Lemma 6.3.} Assume that $u(x)\geq 0$ satisfies that
$-{\rm{div}}(|\nabla u|^{p-2}\nabla u)=\lambda u^{p-1}$ in
$\Omega$ for some $\lambda>0$. Suppose further that for some $t>0$
the level set $\{ x\in\Omega, u(x)>t\}=B_{r(t)}(x_0)$ is a ball
with radius $r(t)$ and center $x_0$. If $u\in
C(\overline{B}_{r(t)}(x_0))$ and $\sigma(\partial_E
B_{r(t)}(x_0))=0$, then $u$ is radially symmetric with respect to
$x_0$ in $B_{r(t)}(x_0)$.}

\vskip 0.2cm

This lemma is crucial to the proof of Theorem 1.2. In the case $p=
2$, the conclusion of the Lemma 6.3 is a famous result due to
Gidas, Ni and Nirenberg \cite{sGNN} (see also Corollary 3.4 in
\cite{sFra}). In the case $1<p<2$, the conclusion of Lemma 6.3 was
given in \cite{sDama}. In the case $p=N$, the conclusion of Lemma
6.3 was proved in \cite{sKeP}.  However, the conclusion of Lemma
6.3 for the case $p>2$ and $p\neq N$ is not available so far.
Here, we give a proof of Lemma 6.3 for all $p\in(1,+\infty)$.

\vskip 0.2cm

{\bf Proof of Lemma 6.3.} By the assumption, we know that for the
same $t>0$ in the above Lemma, $u(x)$ is a solution of the
following Dirichlet problem
\begin{equation}\label{eq47}
\left\{
\begin{array}{ll}
-{\rm{div}}(|\nabla v|^{p-2}\nabla v)=\lambda v^{p-1} &\mbox{in}\ B_{r(t)}(x_0),\\
v\geq 0 & \mbox{in}\ B_{r(t)}(x_0),\\
v=t &\mbox{on}\ \partial B_{r(t)}(x_0).
\end{array}
\right.
\end{equation}
By Lemma 2.4 and Lemma 2.5, we know that any solution of problem
(\ref{eq47}) is strictly positive in $\overline{\Omega}$. Since
problem (\ref{eq47}) is invariant under rotation, we can prove
Lemma 6.3 by proving that uniqueness theorem is valid for
(\ref{eq47}). To this end, we denote $B_{r(t)}(x_0)$ by $\Omega$
for simplicity, and suppose that $v_1>0$ and $v_2>0$ are two
solutions of problem (\ref{eq47}). Then $v_i(x)\ (i=1,2)$ satisfy
\begin{equation}\label{eq48}
\left\{\begin{array}{ll}
-{\rm{div}}(|\nabla v_i|^{p-2}\nabla v_i)=\lambda |v_i|^{p-2}v_i & \mbox{in}\ \Omega,\\
v_i=t & \mbox{on}\ \partial\Omega.\\
\end{array}
\right.
\end{equation}

Let
$$\eta_1=v_1-v_2^pv_1^{1-p}=\frac{v_1^p-v_2^p}{v_1^{p-1}},
\ \ \eta_2=v_2-v_1^pv_2^{1-p}=\frac{v_2^p-v_1^p}{v_2^{p-1}}.$$ It
is obvious that $\eta_i=0\ (i=1,2)$ on $\partial\Omega$.
Multiplying equation (\ref{eq48}) by $\eta_i\ (i=1,2)$ and
integrating by parts, we obtain
$$\int_{\Omega}|\nabla v_i|^{p-2}\nabla v_i\cdot\nabla\eta_i
-\lambda\int_{\Omega}v_i^{p-1}\eta_i=0,\ \ (i=1,2).$$ By a similar
argument to that used in the proof of Proposition 2.2, we infer
that $\nabla (\ln v_1-\ln v_2)=0$, namely, $v_2=Cv_1$ for some
constant $C$. Since $v_1(x)=v_2(x)=t$ for $x\in\partial\Omega$, we
obtain that $C=1$ and $v_1(x)\equiv v_2(x)$ on
$\overline{\Omega}$. Hence, the solution of problem (\ref{eq47})
is unique, and hence, the symmetry result of Lemma 6.3 follows.

\vskip 0.2cm

{\bf Proof of Theorem 1.2.} Let $\Omega$ satisfy
$\lambda_1(\Omega)=\lambda_1(B)$ and $|\Omega|=|B|$, $U_t$ be the
level set of eigenfunction $\psi_{\Omega}$ correspond to
$\lambda_1(\Omega)$. Then by Lemma 6.1,
$H_{\Omega}(U_t,\Phi)=\lambda_1(B)$ for almost all $t\in (m,1)$,
and so,  $U_t$ is a ball for any $t\in(m,1)$ and
$\sigma(\partial_EU_t)=0$ by Lemma 6.2. At this stage, Lemma 6.3
implies that $\psi_{\Omega}$ is radially symmetry inside $U_t$,
and all interior level sets $U_{\tau}$ for $\tau\in(t,1)$ are
concentric balls. In particular, for all $t\in(m,1]$, the level
sets $U_t$ are concentric balls. Therefore, $\Omega=\bigcup_{t\in
(m,1)}U_t$ is a ball.

\vskip 0.5cm

{\bf Acknowledgements}. Heartfelt thanks are given to Professor
Xu-jia Wang for many invaluable comments.


\newpage

\begin{thebibliography}{s2}

\bibitem{sBha} T.Bhattacharya, A proof of the Faber-Krahn inequality
for the first eigenvalue of the p-Laplacian, Annali di Matematica
pura ed applicata, IV 177 (1999), 325-343

\bibitem{sBos} M.H.Bossel, Membranes $\acute{e}$lastiquement
li$\acute{e}$es inhomog$\grave{e}$nes ou sur une surface: une
nouvelle extension du th$\acute{e}$or$\acute{e}$me
isop$\acute{e}$rim$\acute{e}$trique de Rayleigh-Faber-Krahn, Z.
Angew. Math. Phys., 39(1988), 733-742

\bibitem{sBDan}D.Bucur and D.Daners, An alternative approach to the Faber-Krahn inequality
for Robin problems, Calc. Var. PDE., 37(2010), p75-86

\bibitem{sDama}L.Damascelli and F.Pacella, Monotonicity and
symmetry results for p-Laplace equations and applications, Adv.
Differential Equations, 5(2000), 1179-1200

\bibitem{sDan} D.Daners, A Faber-Krahn inequality for Robin problems
in any space dimension, Math. Ann. 335 (2006), 767-785

\bibitem{sDanK} D.Daners, J. Kennedy, Uniqueness in the Faber-Krahn inequality
for Robin problems, SIAM J. Math. Anal., 39(2007), 1191-1207

\bibitem{sDiBen}E.DiBenedetto, $C^{1,\alpha}$ local regularity of
weak solutions of degenerate elliptic equations, Nonlinear Anal.
TMA, 7(1983), 827-850


\bibitem{sFab} G.Faber, Beweis, dass unter allen homogenen Membranen
von gleicher Fl\"{a}che und gleicher Spannung die kreisf\"{o}rmige
den tiefsten Grundton gibt, Sitzungsber, Bayr. Akad. Wiss.
M\"{u}nchen, Math.-Phys. Kl, 1923, 169-127

\bibitem{sFra}L.E.Fraenkel, An Introduction to Maximum Principles
and Symmetry in Elliptic Problems, Cambridge Tracts in Math. 128,
Cambridge University Press, Cambridge, UK, 2000

\bibitem{sGNN} B.Gidas, W.M.Ni and L.Nirenberg, Symmetry and
related properties via the maximum principle, Comm. Math. Phys.,
68(1979), 209-243

\bibitem{sGio} T.Giorgi, R.G.Smits, Monotonicity results for the
principal eigenvalue of the generalized robin problem, Illinois J.
Math. 49 (4) (2005), 1133-1143

\bibitem{sGuoD} D.J.Guo, Nonlinear  Functional Analysis(in Chinese),
Second Editition, Shangdong Science and Technology Press, (2002)

\bibitem{sIssa} Chavel Isaac, Eigenvalues in Riemannian Geometry,
Academic Press, Inc, Orlando, FL, 1984

\bibitem{sKaw} B.Kawohl, Rearrangements and convexity of Level sets
in PDE (Lecture Notes in Mathematics 1150), Springer-Verlag,
Heidelberg 1985

\bibitem{sKeP} S.Kesavan, F.Pacella, Symmetry of positive
solutions of a quasilinear elliptic equation via isoperimetric
inequalities, Applicable Anal., 54(1994), 27-37

\bibitem{sKra} E.Krahn, \"{U}ber eine von Rayleigh formulierte
Minimaleigenschaft des Kreises, Math. Ann. 94 (1925), 97-100

\bibitem{sKra1} E.Krahn, \"{U}ber Minimaleigenschaften der Kugel in drei und mehr
Dimensionen, Acta Comm. Univ. Tartu (Dorpat) A9 (1926), 1-44

\bibitem{sLind} P.Lindqvist, On the equation $-{\rm{div}}(|\nabla u|^{p-2}\nabla u)
=\lambda |u|^{p-2}u $, Proc. Amme. Math. Soci. 109 (1) (1990),
157-163

\bibitem{sSerrin} P.Pucci and J.Serrin, The strong maxium principle
revisited, J. Differential Equations, 196(2004),1-66


\bibitem{sRay} J.W.S. Rayleigh, The Theory of Sound, second edition
revised and enlarged (in 2 volumes). Dover publications, New York,
1945 (republication of the 1894/96 edition)

\bibitem{sSaka} S.Sakaguchi, Concavity properity of solutions to some
degenerate quasilinear elliptic dirichlet problems. Annali della
Scuola Normale Superiore de Pisa, Serie IV (Classe di Scienze) 14,
1987,pp.403-421

\bibitem{sSze} G.Szeg\"{o}, Inequalities for certain eigenvalues of a
membrane of given area, J. Rational Mech. Anal. 3 (1954), 343-356

\bibitem{sTolk}P.Tolksdorf, Regularity for a more general class of
quasilinear elliptic equations, J. Differential Equations,
51(1984), 126-150

\bibitem{sTolk1}P.Tolksdorf, On the Dirichlet problem for
quasilinear equations in domains with conical boundary points,
Comm. Partial Differential Equations, 8(1983), 773-817


\bibitem{sWein} H.F.Weinberger, An isoperimetric inequality for the
$n$-dimensional free membrane problem, J. Rational Mech. Anal. 5
(1956), 633-636

















\end{thebibliography}
\end{document}